\documentclass[12pt]{amsart}
\usepackage{amsfonts}
\usepackage{amssymb}

\theoremstyle{definition}

\theoremstyle{remark}

\begin{document}

\title[ On the holonomy group of hypersurfaces ]
{On the holonomy group of hypersurfaces of spaces of constant curvature}%

\author{Ognian Kassabov}%
\address{Department of Mathematics and Informatics}
\address{University of Transport}
\address{158 Geo Milev Str.}
\address{1574 Sofia, Bulgaria}

\email{okassabov@abv.bg}


\begin{abstract}
We classify hypersurfaces $M^n$ of manifolds of constant nonzero sectional 
curvature  according their restricted homogeneous holonomy groups.
It turns out that outside of the evident cases (restricted holonomy 
group  $SO(n)$ and flat submanifolds) only two cases arise: restricted holonomy group  
$  {SO}(k)\times  {SO}(n-k)$  (when  $ M $
is locally a product of two space forms) and $  {SO}(n-1) $ (when 
$ M $ is locally a product of an
$ (n-1)$-dimensional space form and a
segment).
\end{abstract}

\maketitle

\section{Introduction}

The holonomy groups are fundamental analytical objects
in the theory of manifolds and especially in the theory of Riemannian manifolds. 
The holonomy group of a Riemannian manifold reflects for example on local 
reducibility of the manifold. In \cite{Kurita} M. Kurita classifies the conformal flat 
Riemannian manifolds according their restricted homogeneous holonomy group. 

There exists a similarity between the conformal flat Riemannian manifolds and the
hypersurfaces of a Riemannian manifold, see e.g. a remark of R. S.
Kulkarni in \cite{Kulkarni}. So it is natural to look for a result in the submanifold geometry,
analogous to the Kurita's theorem.
In \cite{Kobayashi} S. Kobayashi proves that the holonomy group
of a compact hypersurface of $\bf E$$^{n+1}$ is $ {SO}(n)$.
Generalizations of of Kobayashi's result are obtained by R. Bishop \cite{Bishop} and G. Vranceanu \cite{Vranceanu}. 

In this paper we consider analogous question for hypersurfaces of non-flat real space forms
according their holonomy groups. Namely we prove:

\vspace{2 mm}
{\bf Theorem 1.} {\it Let $ M^n $ ($n\ge 3$)  be a connected hypersurface of a
space $ \widetilde M^{n+1}(\nu )$ of constant positive sectional curvature $ \nu $.
Then the restricted homoge\-neous holonomy group $ H_p $ of $ M^n $ in
any point $ p $ is in general the special orthogonal group $ {SO}(n) $.
If  $ H_p  $ is not $ SO(n) $ at any point $p \in M^n$, then one of
the following cases appears:

a) $ H_p =  {SO}(k)\times  {SO}(n-k)$, $ 1 < k < n-1 $ and $ M^n $
is locally a product of a $ k $-dimensional space of constant curvature
$ \nu + \lambda^2 $ and an $ (n-k) $-dimensional space of constant
sectional curvature $ \nu + \mu^2 $, with $ \nu + \lambda \mu = 0 $;

b) $ H_p = {SO}(n-1) $ and $ M^n $ is locally a product of an
$ (n-1)$-dimensional space of constant sectional curvature and a
segment.}

\vspace{2 mm}
A similar theorem for complex manifolds is proved in \cite{N-S}.

\section{ Preliminaries.}
\noindent\par
Let $ \widetilde M^{n+1} $ be an $ (n+1) $-dimensional Riemannian manifold
with metric tensor $ g $ and denote by $ \widetilde\nabla $ its
Riemannian connection. It is well known that if $\widetilde M^{n+1} $ is of
constant sectional curvature $ \nu $, then its curvature operator
$ \widetilde R $ has the form
$$
	\widetilde R(x,y) = \nu\, x\land y     \ \ \ ,
$$
where the operator $ \land $ is defined by
$$
	(x\land y)z = g(y,z)x - g(x,z)y \ \ \  .
$$
Such a manifold is denoted by $ \widetilde M^{n+1}  (\nu) $. Now let $ M^n $
be a hypersurface of $ \widetilde M^{n+1}(\nu) $ and denote by $ \nabla $ its
Riemannian connection. Then we have the Gauss formula
$$
	\widetilde \nabla _X Y = \nabla _X Y + \sigma (X,Y)
$$
for vector fields $ X, Y $ on $ M^n $, where $ \sigma $ is a
normal-bundle-valued symmetric tensor field on $ M^n $, called the
second fundamental form of $ M^n $ in $ \widetilde M^{n+1} $. Let $ \xi $
be a unit normal vector field. Then the Weingarten formula is
$$
	\widetilde\nabla_x\xi = - A_{\xi}X
$$
and the operator $ A_{\xi} $ is related to $ \sigma $ by
$$
   g(\sigma (X,Y),\xi) = g(A_{\xi}X,Y)=g(A_{\xi}Y,X)    \ \ \ .
$$
Suppose that we have fixed a normal vector field $ \xi$. Then we shall
write $ A $ insteed of $ A_{\xi} $. The equations of Gauss and
Codazzi are given respectively by
$$
	R(X,Y) = \nu (X\wedge Y) + AX\wedge AY  \ \ \ ,
$$
$$
	(\nabla_XA)Y = (\nabla_YA)X        \ \ \ ,
$$
$ R $ denoting the curvature operator of $ M^n $.

It is known that the Lie algebra of the infinitesimal holonomy group at a point $p$ of
a Riemannian manifold $M$ is generated by all endomorphisms of the form
$$
	(\nabla^kR)(X,Y;V_1,...,V_k) \ ,
$$ 
where $X,Y,V_1,...,V_k \in T_pM$ and $0\le k<+\infty$ \cite{K-N}. Moreover if the
dimension of the infinitesimal holonomy group is constant, this group
coincides with the restricted homogeneous holonomy group \cite{K-N}.

\section{ Proof of  Theorem 1.}
\noindent\par
Let $ p $ be an arbitrary point of $ M^n$. We choose an orthonormal basis
$ {e_1,...,e_n} $ of $ T_pM $, which diagonalize the symmetric
operator $ A $, i.e.
$$
	Ae_i = \lambda_ie_i \ \ \ \ \ \ i=1,...,n \ .
$$
Then by the equation of Gauss we obtain
$$
	R(e_i,e_j) = (\nu + \lambda_i\lambda_j)e_i\land e_j \ \ \ . \leqno (3.1) 
$$

First we note that $ M^n $ cannot be flat at $ p $. Indeed if $ M^n $
is flat, we obtain  from (3.1) \ $ \nu + \lambda_i\lambda_j = 0 $ for all $ i \ne j $.
Since $n>2$ this implies easily $ \nu + \lambda_1^2 = 0 $, and because of $\nu>0$ this is a
contradiction.

Since $ M^n $ is not flat at $ p $, there exist $ i \ne j $, such
that $\nu + \lambda_i\lambda_j \ne 0 $. Then (3.1) implies that
$ e_i \land e_j $ belongs to the Lie algebra $ h_p $ of $ H_p $.
As in \cite{Kurita} we denote by $ {SO} [i_1,...,i_k] $ the subgroup of
$  {SO} (n) $, which induces the full rotation of the linear
subspace, generated by $ e_{i_1},...,e_{i_k} $ and fixes the
remaining vectors. Denote also by  $  {so} [i_1,...,i_k] $
the Lie algebra of $  {SO} [i_1,...,i_k] $. Then according to
the above argument $ H_p $ contains  $  {SO} [i,j] $.

If $ H_p$ contains $ {SO} (n) $, then $ H_p =  {SO} (n) $,
because the restricted homogeneous holonomy group  $ H_p$  of a 
Riemannian manifold is a subgroup of $ {SO} (n), $ see  \cite{Bor-Lichn}.

Let $ H_p $ is not $ {SO} (n) $. Then there exist $ k $,
$ 2 \le k \le n-1 $ and indices  $ i_1,...,i_k $, such that
$ H_p $ contains  $  {SO} [i_1,...,i_k] $ but doesn't
contain $  {SO} [i_1,...,i_k,u ] $ for $ u \ne i_1,...,i_k $.
Without loss of generality we can assume that $ H_p $ contains $ SO [1,...,k] $,
but does not contain $ SO [1,...,k,u ] $ for $ u > k $.

Let us suppose that $ h_p $ contains $  so [a ,u ] $ for some
$ a \in \{ 1,...,k\} $ and $ u \in \{ k+1,...,n \} $. Since
$$
	[e_b \wedge e_a ,e_a \wedge e_u] = e_b \land e_u
$$
it follows that the Lie algebra $ h_p $ contains $ e_b \wedge e_u $ for
$ b = 1,...,k $. Hence $ h_p $ contains $  so [1,...,k,u ] $,
which is a contradiction.

Consequently $ h_p $ doesn't contain  $  so [a,u] $ for any
$ a = 1,...,k $;   $ u =k+1,...,n $. Then (3.1) implies
$$
	\nu + \lambda _a \lambda _u = 0 \ \ \ \ a=1,...,k;\
	u = k+1,...,n .    \leqno (3.2)
$$
Hence, using $\nu \ne 0$, we obtain $ \lambda_1 = ... = \lambda_k $ and
$ \lambda _{k+1} =...= \lambda _n $. Denote $ \lambda = \lambda_1 $;
$ \theta = \lambda _{k+1} $. Then by (3.2) \ $ \nu + \lambda\theta = 0 $, $\lambda\ne 0$, $\theta\ne 0$
and it follows easily $ \lambda \ne \theta $, $ \nu + \lambda^2 \ne 0 $, 
$ \nu + \theta^2 \ne 0 $. 

In a neighborhood $ W $ of $ p $ we consider continuous functions
$ \Lambda_1,...,\Lambda_n $, such that for any point $ q \in W $
the numbers $ \Lambda_1 (q),...,\Lambda_n (q) $ are the eigenvalues
of $ A $. Since $ \nu + \lambda^2 \ne 0 $, $ \nu + \theta^2 \ne 0 $,
then in an open subset $ V $ of $ W $ containing $ p $ we have
$$
	\nu + \Lambda_a(q)\Lambda_b(q) \ne 0 \ \ \ \  a,b=1,...,k \ \ \ ;
$$
$$    \ \ \
	\nu + \Lambda_u (q)\Lambda_v (q) \ne 0 \ \ \ \
	u , v = k+1,...,n   \ \ \ .
$$
Hence $ H_q $ contains $  SO [1,...,k] $ and
$ SO [k+1,...,n] $. Suppose that $\nu+\Lambda_a(q)\Lambda_u(q)\ne 0$
for some $a=1,...,k$, $u=k+1,...,n$. Then $h_q$ contains
$e_a \wedge e_u$, so as before $h_q$ contains $so[1,...,k,u]$ and
analogously $h_q$ contains $so(n)$, which is not possible. So
$ \nu + \Lambda_a(q)\Lambda_{\alpha }(q) = 0 $. Hence as before
we find
$$
	\Lambda_1(q)=...=\Lambda_k(q) \ \  ,  \ \ \
	\Lambda_{k+1}(q)=...=\Lambda_n(q) \ \ \ .
$$
\par
Consequently in a neighborhood $ V $ of $ p $ there exist a number
$ k $ and continuous functions $ \Lambda(q) , \Theta(q) $ such that
$ \Lambda(q) \ne \Theta(q)$ and
$$
	\Lambda_1(q)=...=\Lambda_k(q)= \Lambda (q)\ne 0 \  \ , \ \ \
	\Lambda_{k+1}(q)=...=\Lambda_n(q)=\Theta (q)\ne 0 \ \ \       \leqno (3.3)
$$
for $ q \in V $. Since $M^n$ is connected \ $k$ \ is a constant on $ M^n $.
Consequently (3.3) holds on $ M^n $. On the other hand using $ \nu + \Lambda\Theta = 0 $
and the fact that \ $ k\Lambda + (n-k)\Theta = tr A $ \ is smooth we conclude that
$ \Lambda $ and $ \Theta $ are smooth functions on $ M^n $. Define
two distributions
$$
	T_1(q) = \{ x \in T_q(M) \ : \ \ Ax = \Lambda (q)x \}  \ \ \ ,
$$
$$
	T_2(q) = \{ x \in T_q(M) \ : \ \ Ax = \Theta (q)x \}   \ \ \ .
$$
It follows directly that $ T_1 $ and $ T_2 $ are orthogonal and for $ X,Y \in T_1 $,
$ Z,U \in T_2 $ we have
$$
	R(X,Y) = (\nu + \Lambda^2)X\wedge Y \ \ \ ,
$$
$$
	R(Z,U) = \frac{\nu}{\Lambda^2}(\nu + \Lambda^2)Z\wedge U \ \ \ ,
$$
$$
	R(X,Z) = 0  \ \ \ .
$$
We choose local orthonormal frame fields $ \{ E_1,...,E_k \} $  of $ T_1 $
and $ \{ E_{k+1},...,E_n \} $  of $ T_2 $ and we denote
$$
	\nabla _{E_i}E_j = \sum _{s=1}^{n} \Gamma_{ijs}E_s   \ \ \ .
$$
Then $ \Gamma_{ijs} = -\Gamma_{isj} $ for all $ i,j,s = 1,...,n $,
in particular $ \Gamma_{ijj} = 0 $. As before let
$ a,b,c \in \{ 1,...,k\} $ and $ u,v \in \{ k+1,...,n\} $. 
From the second Bianchi identity we have
$$
	(\nabla_aR)(E_b,E_u) +(\nabla_bR)(E_u,E_a)+(\nabla_u R)(E_a,E_b) = 0
$$
and hence
$$   
		\begin{array} {rl}   \displaystyle
        E_u(\Lambda^2) E_a\land E_b 
        + (\nu + \Lambda^2)  \sum _{c=1}^{k}\left\{ 
                         \Gamma_{buc} E_a \land E_c - \Gamma_{auc} E_b \land E_c \right\}  & \\
        \displaystyle +(\nu+\Lambda^2) \sum_{v =k+1}^{n} \left\{\frac {\nu }{\Lambda^2}
           (\Gamma_{abv} - \Gamma_{bav}) E_u\land E_v +
        \Gamma_{uav} E_v\land E_b - \Gamma_{ubv} E_v\land E_a  \right\} & = 0     \ \ \ .
    \end{array} 
$$
Consequently we obtain
$$
	E_u(\Lambda^2) = (\nu + \Lambda^2)
	\{\Gamma_{aau}+\Gamma_{bbu} \}   \ \ \ , \leqno (3.4)
$$
$$
	(\nu + \Lambda^2)\Gamma_{uva} = 0
$$
for all $ a \ne b $. Since $ \nu + \Lambda^2 \ne 0 $ we find
$ \Gamma_{uva} = 0 $, so $ T_2 $ is parallel.

Let $n-k\ge 2$. Then analogously to the above $ T_1 $ is also parallel.
Now (3.4) implies that $ \Lambda $ doesn't depend on $ E_u $
and analogously $ \Theta  $ doesn't depend on $ E_a $. Hence, using
$\nu+\Lambda\Theta=0$ we conclude that $\Lambda $ and $ \Theta $ are constants. 
So we obtain the case a) of our Theorem.

Let $ n-k =1 $. We shall show that
under the assumption $ H_p \ne  SO (n) $ the distribution
$ T_1 $ is again  parallel. By the Codazzi equation we have
$$
    (\nabla_aA)(E_b) =(\nabla_bA)(E_a)   \ \ \ .
$$
This implies
$$
	E_a(\Lambda )E_b + (\Lambda - \Theta)\Gamma_{abn}E_n =
	E_b(\Lambda )E_a + (\Lambda - \Theta)\Gamma_{ban}E_n    \ \ \ .
$$
Hence $ E_a(\Lambda ) = 0 $ for $a=1,...n-1$. Now from
$$
	(\nabla_aA)(E_n) =
	(\nabla_nA)(E_a)   \ \ \
$$
we obtain
$$
	E_n(\Lambda ) E_a + (\Lambda - \Theta )\sum_{c=1}^{n-1}
\Gamma_{anc}E_c = 0    \ \ \ .
$$
Hence we derive
$$
  E_n (\Lambda ) = (\Lambda - \Theta )\Gamma_{aan} \ , \leqno (3.5)
$$
$$
  (\Lambda - \Theta )\Gamma_{acn}=0  \qquad {\rm for}\ c\ne a\ .
$$
Since $\Lambda\ne \Theta$ the last equality implies $ \Gamma_{acn} = 0 
$ for $a \ne c $. On the other hand (3.5)
implies $ \Gamma_{aan} = \Gamma_{bbn} $. If $ \Gamma_{aan} = 0 $,
then $ T_1 $ is parallel and from (3.5) $E_n(\Lambda)=0$, so $\Lambda$ is 
a constant. Because of $\nu+\Lambda\Theta\ne 0$ it follows that $\Theta$ is
a constant too. Hence we obtain the case b) of our
Theorem. Let us suppose that $ \Gamma_{aan} \ne 0 $. We compute
directly
$$
	(\nabla_aR)(E_a,E_b) = (\nu + \Lambda^2)\Gamma_{aan}E_n \land E_b \ \ \ .
$$
Hence $ E_n \land E_b \in h_p $ and as before it follows that
$ SO (n) = H_p $, which is not our case. This proves Theorem 1.

\vspace{2mm}
{\bf Remark.} In the same way we can consider the case where
$ \widetilde M^{n+1}(\nu ) $ is of constant negative sectional curvature
$ \nu $. Then we obtain

\vspace{2 mm}
{\bf Theorem 2.} {\it Let $ M^n $ ($n\ge 3$)  be a connected hypersurface of a
space $ \widetilde M^{n+1}(\nu )$ of constant negative sectional curvature $ \nu $.
Then the restricted homoge\-neous holonomy group $ H_p $ of $ M^n $ in
any point $ p $ is in general the special orthogonal group $ {SO}(n) $.
If $M^n$ is not flat and  $ H_p $ is not $ SO(n) $ at any point $p \in M^n$, 
then one of the following cases appears:

a) $ H_p =  {SO}(k)\times  {SO}(n-k)$, $ 1 < k < n-1 $ and $ M $
is locally a product of a $ k $-dimensional space of constant curvature
$ \nu + \lambda^2 $ and an $ (n-k) $-dimensional space of constant
sectional curvature $ \nu + \mu^2 $, with $ \nu + \lambda \mu = 0 $

b) $ H_p = {SO}(n-1) $ and $ M $ is locally a product of an
$ (n-1)$-dimensional space of constant sectional curvature and a
segment.}

\end{document}